\documentclass[a4, amsmath, amsthm, amssymb]{amsart}

\usepackage{pdfpages}
\usepackage[left=20mm,top=0.6in,bottom=12mm]{geometry}

\numberwithin{equation}{section}

\newtheorem{Lemma}{LEMMA}[section]
\newtheorem{Theorem}[Lemma]{Theorem}
\newtheorem{Proposition}[Lemma]{Proposition}
\newtheorem{Corollary}[Lemma]{Corollary}

\newtheorem{remark}[Lemma]{Remark}
\newtheorem{definition}[Lemma]{Definition}
\newtheorem{example}[Lemma]{Example}

\newtheorem{Fact}[Lemma]{Fact}
\newtheorem{assumption}[Lemma]{Assumption}

\def\bt{\begin{Theorem}}
\def\et{\end{Theorem}}
\def\bl{\begin{Lemma}}
\def\el{\end{Lemma}}
\def\bp{\begin{Proposition}}
\def\ep{\end{Proposition}}
\def\bcor{\begin{Corollary}}
\def\ecor{\end{Corollary}}
\def\bpf{\begin{proof}}
\def\epf{\end{proof}}

\def\brem{\begin{remark}\rm }
\def\erem{\hfill $\lozenge$ \end{remark}}

\def\bedef{\begin{definition}\rm }
\def\endef{\hfill$\lozenge$\end{definition}}

\def\beg{\begin{example}\rm }
\def\eeg{\hfill $\lozenge$\end{example}}

\def\bef{\begin{Fact}}
\def\eef{\end{Fact}}

\def\bea{\begin{assumption}}
\def\ena{\end{assumption}}
\def\bc{\begin{center}}
\def\ec{\end{center}}
\def\noi{\noindent}

\def\beq{\begin{equation}}
\def\eeq{\end{equation}}
\def\beqarray{\begin{eqnarray*}}
\def\eeqarray{\end{eqnarray*}}
\def\<{\leftangle}
\def\>{\rightangle}
\def\({\left(}
\def\){\right)}
\def\f{\varphi}

\def\<{\langle}
\def\>{\rangle}
\def\q{\quad}

\def\a{\alpha}

\def\g{\gamma}

\def\d{\delta}

\def\t{\tau}

\def\e{\varepsilon}

\def\O{\Omega}

\def\w.r.t.{with respect to}
\def\R{{\mathbb{R}}}

\def\bq{\begin{quote}}
\def\eq{\end{quote}}

\def\bit{\begin{itemize}}
\def\eit{\end{itemize}}

\def\ben{\begin{enumerate}}
\def\een{\end{enumerate}}

\begin{document}
\title[identification of impedance coefficient]{Convergence rates for identification of Robin coefficient from terminal observations}
\author{Subhankar Mondal }
\address{TIFR Centre for Applicable Mathematics, Bangalore-560065, India}
\email{
subhankar22@tifrbng.res.in; s.subhankar80@gmail.com
}
\maketitle

\begin{abstract} This paper deals with the problem of identification of a Robin coefficient (also known as impedance coefficient) in a parabolic PDE from terminal observations of the temperature distributions. The problem is ill-posed in the sense that small perturbation in the observation may lead to a large deviation in the solution. Thus, in order to obtain stable approximations, we employ the Tikhonov-regularization. We propose a weak source condition motivated by the work of Engl and Zou (2000) and  obtain a convergence rate of $O(\d^\frac{1}{2})$, the main goal of this paper, where $\d$ is the noise level of the observed data. The obtained rate is better than some of the previous known rates.
Moreover, the advantage of the proposed source condition is that we are getting the above mentioned convergence rate without the need for characterizing the range space of modelling operator, which is in contrast to the general convergence theory of Tikhonov-regularization for non linear operators, where one obtain the same order of convergence by characterizing the range of the adjoint of the Fr{\'e}chet derivative of modelling operator, a challenging task for many problems.    
\end{abstract}

\section{Introduction}

Let $d\in \{2,3\}$ and $\O\subset \R^d$ be a bounded domain with Lipschitz boundary $\partial\O$.
Let $\tau>0$ be fixed and we denote the sets $\O\times [0,\tau]$ and $\partial\O\times [0,\tau]$ by $\O_\tau$ and $\partial\O_\tau$, respectively.
We consider the PDE 
\beq\label{pde1}
 \begin{cases}
    u_t-\Delta u=f\q &\text{in}\q \O_\tau,\\
    \frac{\partial u}{\partial \nu }+\gamma(x)u= g\q &\text{on}\q \partial\O_\tau,\\
    u(\cdot, 0)=u_0\q &\text{in}\q \O,
 \end{cases}
\eeq
where $f\in L^2(0,\tau;L^2(\O)),\, \,g\in L^2(0,\tau;L^2(\partial\O)),\,\,\g\in L^\infty(\partial\O)$ and $u_0\in L^2(\O).$ Here, for a Banach space $Y$, we used the notation  $L^2(0, \tau; Y)$ for the space of all $Y$-valued measurable functions $\phi$ on $[0, \tau]$ such that $\int_0^\tau \|\phi(t)\|_Y^2 dt<\infty$. The Robin coefficient $\g$ is the impedance coefficient that represents the heat exchange on the boundary $\partial\O,\,\,\nu$ denotes the outward unit normal of the boundary $\partial\O.$ The direct or forward problem for \eqref{pde1} is to find $u(x,t)$ satisfying \eqref{pde1}(possibly in weak sense) for the known impedance coefficient $\g$ and the input data $f,g$ and $u_0$, and the existence of solution of the forward problem is well known. 

The system \eqref{pde1} models heat conduction phenomenon where the impedance coefficient characterizes the thermal properties of the conductive material on the interface and certain physical processes, e.g. corrosion, on the boundary \cite{bellassoued_cheng_choulli_2008, jin_zou_2009}. Thus, the value of the impedance coefficient $\g$ is of significant interest in thermal imaging such as safety analysis of nuclear reactor and thermal protection of space shuttles \cite{beck_blackwell_clair}. In practice, the impedance coefficient $\g$ cannot be specified from direct measurements since the domain $\O$ may be embedded in an unknown region \cite{wang_liu_2017}. Therefore, one has to deal with the inverse problem of identifying the impedance coefficient from some available observations, for example, a final time observation on the whole spatial domain \cite{wang_liu_2017}, partial Dirichlet boundary observation for the full time period \cite{bellassoued_cheng_choulli_2008, jin_lu_2012}, time integral observation on the full spatial boundary \cite{hao_thanh_lesnic_2013, liu_wang_2016}, final or an intermediate time observation on the full spatial boundary \cite{hao_thanh_lesnic_2013}. Considering the amount of work that has been devoted for this type of parameter identification problem and various type of observations over the years, it is impossible to list all of them, however, the interested reader may refer to \cite{cheng_lu_yamamoto_2012, isakov_1991, jin_zou_2009, keung_zou_1998, kugler_sincich_2009, yamamoto_zou_2001, zhang_liu_2021}.

In the applications of heat transport in the high temperature, it is not possible to measure  the temperature distribution in the whole time interval $[0,\t].$ Therefore, following \cite{engl_zou_2000, cao_pereverzev_2006} (see also the recent work \cite{k.cao_2022}), in this paper, we assume that the terminal status observations of the temperature distribution is known and with this knowledge we consider the inverse problem of identifying the spatially dependent impedance coefficient $\g$. More precisely, we assume that 
\beq\label{exact_obs_data}
u(x,t)=\phi(x,t)\q\text{in}\,\,\O\times[\tau-\sigma, \tau],
\eeq
is known at hand, where $\sigma>0$ is small so that $\t-\sigma>0,$ and then we consider the problem of identifying $\g$ from the knowledge of $\phi$ or its noisy approximations $\phi_\d$ satisfying
 \beq\label{noisy_obs_error}
 \int_{\tau-\sigma}^\tau \|\phi-\phi_\delta\|^2_{L^2(\O)}\,dt\leq \delta^2
 \eeq
for some noise level $\d>0.$

It can be observed that our inverse problem is non linear. Furthermore, as can be seen from next section, the inverse problem is also ill-posed, that is, a small perturbation in the observed data may lead to a large deviation in the corresponding solutions. Thus, some regularization scheme has to be employed in order to obtain stable approximations. We employ the standard Tikhonov regularization for obtaining the approximations for $\g$. It is well known that the convergence of the stable approximations obtained by regularization can be arbitrary slow \cite{schock_1984} unless some apriori conditions, the so-called source conditions (cf. \cite{engl_hanke_neubauer, nairopeq}) is assumed on the unknown that has to be identified. In Tikhonov regularization theory for non linear operators in Hilbert spaces, generally the source condition involves the adjoint of the Fr{\'e}chet derivative of the non linear operator involved \cite{engl_hanke_neubauer}. More precisely, if $\mathcal{X}$ and $\mathcal{Y}$ are Hilbert spaces, $F:\mathcal{X}\to \mathcal{Y}$ is a (non linear) operator which is Fr{\'e}chet differentiable, consider the problem of solving an ill-posed operator equation $$F(x)=y$$ in the sense that small perturbation in $y$ may lead to a large deviation in the solution of the operator equation. Let $y$ be the exact data and $x^\dagger$ be the unique solution (i.e. $F(x^\dagger)=y$) to be identified from the knowledge of $y^\text{obs}\in \mathcal{Y}$, the observed data satisfying $\|y-y^\text{obs}\|_\mathcal{Y}\leq \d.$ The stable approximations are the minimizers of the functional 
$$\min_{x\in\mathcal{X}}\|F(x)-y^\text{obs}\|^2_{\mathcal{Y}}+\a\|x-x^*\|^2_{\mathcal{X}},$$ for a fixed $\a>0$, the regularization parameter, and $x^*$ is an initial guess for $x^\dagger$ that incorporates some apriori smoothness assumptions. Let $F'$ denotes the Fr{\'e}chet derivative of $F$. Moreover, we assume that the Fr{\'e}chet derivative of $F$ is Lipschitz continuous with the Lipschitz constant $C_{\rm Lip}$. Then it is known that(cf. \cite{engl_hanke_neubauer, engl_kunisch_neubauer_1989, neubauer_1989}) if $$x^\dagger-x^*=F'(x^\dagger)^*\f$$ for some $\f$ satisfying the smallness condition
\beq\label{source_smallness}
C_{\rm Lip}\|\f\|\leq 1
\eeq
then the rate of convergence of the regularized solution is $O(\d^{\frac{1}{2}})$, provided the regularization parameter $\a$ is chosen as $\a\sim\d$. In many problems it is extremely difficult to characterize these range spaces whereas in many problems these range spaces turns out to be certain Sobolev spaces with higher smoothness, see for e.g. \cite{hohage_1997}. In addition to these difficulties, verification of the smallness condition \eqref{source_smallness} is another challenging task. In order to overcome these challenges, a new type of source condition was proposed by Engl and Zou in \cite{engl_zou_2000} for a parameter identification problem in heat conduction, which is simple and verifiable (atleast for some reasonable regularity assumptions on $f,\,g$ and $u_0$), does not require any smallness condition and also does not require much higher smoothness assumption on the unknowns. Motivated by the work in \cite{engl_zou_2000} and also the recent work in \cite{k.cao_2022}, we consider a similar weak source condition and obtain the convergence rates, which is the main goal of this paper. Also, we show explicitly that our source condition is verifiable under certain regularity assumption(see Theorem \ref{source_verify_construct}).

We now discuss the advantages and shortcomings associated with our considered observation \eqref{exact_obs_data} in comparison to some of the recent works where the observations are different from \eqref{exact_obs_data}. 

In \cite{liu_wang_2016} the authors have considered the problem of reconstruction of $\g(x)$ associated with a system similar to \eqref{pde1} with $f=0,\,u_0=0$ from the non-local measurement of the form 
\beq\label{non-local-bdry_obs}
\int_0^\tau w(t)u(x,t)\,dt=h(x)\q\text{on}\,\,\partial\O,
\eeq
for some weight function $w.$ Since the governing PDE is a homogeneous heat equation, by using the fundamental solution of the heat equation, the authors could make use of their boundary observation \eqref{non-local-bdry_obs} in the analysis of the inverse problem of reconstructing $\g$. In fact, the inverse problem of reconstructing $\g$ is transformed into a problem of solving a system of ill-posed non linear integral equations, where at first one has to solve for a certain potential $q(x,t)$ that arises from the fundamental solution (see \cite[pg.4]{liu_wang_2016}) and then solve for the impedance coefficient. Moreover, for the stable reconstruction the authors have considered a semi-Tikhonov regularization scheme in which the penalty term involves the potential $q$ (an auxiliary unknown) in contrast to our approach of traditional Tikhonov regularization functional where the penalty term is comprised of $\g$, the actual unknown to be identified. Although the observation \eqref{non-local-bdry_obs} seems to be more realistic (as it deals with measurement only in the spatial boundary) than \eqref{exact_obs_data}, but their analysis requires $f=0$ and $u_0=0$, which is not the case in this paper. Moreover, the work in \cite{liu_wang_2016} does not provide any error estimate.

In \cite{wang_liu_2017} the authors have considered the the problem of simultaneous identification of $\g(x)$ and the initial temperature $u_0(x)$ associated with a system similar to \eqref{pde1} from the final  time observation on the full spatial domain, that is, the data used for the inversion is of the form
\beq\label{final_time_obs_full_spatial_domain}
u(x,\tau)=h(x)\q\text{in}\,\,\O,
\eeq
and in this work the authors could overcome the restriction of $u_0=0$ that they have considered in their earlier work \cite{liu_wang_2016}. Because of the ill-posedness of the inverse problem, the authors in \cite{wang_liu_2017} have obtained a stable approximations for $\g$ by a regularization scheme that involves the mollification of  inversion input data $h$, and that is achieved by using higher regularity on $h$, namely $h\in W^{3,p}(\O)$ for some $p>2.$ Since $h$ is an inversion input data, assuming such higher regularity of $h$ is not that much realistic from application point of view. Moreover, for a noisy observation $h^\d$ of $h$ the authors have obtained a H{\"o}lder rate of convergence $O(\d^\nu)$, for some $\nu \leq \frac{1}{5}$(see \cite[pg.603]{wang_liu_2017}) under some source condition on $\g$ and $u_0$ which may be not feasible in applications, because it is assumed that the unknowns $\g$ and $u_0$ are sufficiently regular so that $h\in W^{3,p}(\O).$ In contrast to these, as mentioned earlier, in this paper we will obtain stable approximations for $\g$ using Tikhonov regularization and our regularity assumption on the inversion input data $\phi$(see \eqref{exact_obs_data}) is only that $\phi\in L^2(\O\times [\tau-\sigma,\tau]).$ Moreover, under a verifiable source condition that only requires $\g\in H^\frac{1}{2}(\partial\O)$ we obtain a better rate of convergence, namely, $O(\d^{\frac{1}{2}})$. Thus, again it is to be noted that although the observation \eqref{final_time_obs_full_spatial_domain} apparently seems more realistic than \eqref{exact_obs_data}, but in terms of source condition and regularity assumption our work is more realistic than \cite{wang_liu_2017} with a better convergence rate. 

This paper is organised as follows: In Section \ref{sec-2} we collect all the existing results related to existence and uniqueness of solutions for forward problem, precisely state the inverse problem that we consider and analyze the existence, uniqueness and ill-posedness of the inverse problem. In Section \ref{sec-convergence} we do the convergence analysis of the regularized approximations, propose the source condition and proof the convergence rate result, the main result of this paper. In Section \ref{sec-source_diss} we discuss about the source condition, its compatibility, regularity and then show that it is indeed verifiable by a simple construction.

\section{The Inverse Problem}\label{sec-2}
In this section we recall all the definitions, results related to the PDE \eqref{pde1} that will be used later and in addition we formulate the inverse problem more precisely and discuss about its uniqueness and ill-posedness. Throughout the paper whenever we come across a function defined on the boundary $\partial\O$, it is to be understood in the sense of trace \cite{adams, evans}. 
\bedef\label{weak_sol}{\bf (Weak solution)}
An element $u\in L^2(0,\tau;H^1(\O))\cap L^\infty(0,\t;L^2(\O))$ is said to be a weak solution of \eqref{pde1} if
$$
\begin{cases}
&\int_0^\tau \int_\O [-u\frac{\partial\eta}{\partial t}+\nabla u\cdot\nabla \eta]\,dx\,dt+\int_0^\tau\int_{\partial\O}\g u\eta\,dx\,dt\\
&\q=\int_0^\tau \int_\O f\eta\,dx\,dt+\int_0^\tau\int_{\partial\O}g\eta\,dx\,dt+\int_\O\,u_0\eta(\cdot,0)\,dx
\end{cases}
$$
for all $\eta\in H^1(0,\tau;H^1(\O))$ with $\eta(\cdot,\tau)=0.$
\endef

We now state a result about existence and uniqueness of the forward problem. 
Let $$\mathcal{A}=\{\g\in L^\infty(\partial\O):0<\underline{\g}\leq \g\leq \overline{\g}\}$$ be the set of admissible parameters, for some constants $\underline{\g}$ and $\overline{\g}$.
\bt\label{existence_estimate_weak_sol}{\rm (cf. \cite{troltzsch})}
Let $\g\in \mathcal{A}$, $f\in L^2(0,\t;L^2(\O))$, $g\in L^2(0,\tau;L^2(\partial\O))$ and $u_0\in L^2(\O).$ Then there exists a unique weak solution $u$ of \eqref{pde1} satisfying the estimate
\beq\label{weak_sol_estimate}
\max_{[0,\tau]}\|u(t)\|_{L^2(\O)}+\|u\|_{L^2(0,\tau;H^1(\O))}\leq C(\|f\|_{L^2(0,\tau;L^2(\O))}+\|g\|_{L^2(0,\t;L^2(\partial\O))}+\|u_0\|_{L^2(\O)});
\eeq
 where $C$ is a constant depending only on $\underline{\g},\, \O$ and $\t.$
\et
Before proceeding further, let us first precisely state the inverse problem that is considered. 
\begin{quote}
\noi
{\bf (IP)} Identify $\g\in \mathcal{A}$ from the observation $\phi\in L^2(\t-\sigma,\t;L^2(\O))$ such that the unique weak solution $u$ of \eqref{pde1} satisfies \eqref{exact_obs_data}.
\end{quote}
We now discuss about the existence and uniqueness of the solution of the inverse problem {\bf (IP)}.
Let $g\in L^2(\t-\sigma,\t;L^\infty(\partial\O))$, $\phi\in L^2(\t-\sigma,\t;H^1(\O))$ be such that $\phi, \frac{\partial\phi}{\partial\nu}\in L^2(\t-\sigma,\t; L^\infty(\partial\O))$ and $$\int_{\t-\sigma}^\t \phi\,dt\neq 0\q\text{a.e. on}\,\,\partial\O.$$ Then from the Robin boundary condition in \eqref{pde1} and using \eqref{exact_obs_data}, we have
\beq\label{gamma_rep}
\g = \frac{\int_{\t-\sigma}^\t\,g\,dt-\int_{\t-\sigma}^\t\frac{\partial \phi}{\partial\nu}\,dt}{\int_{\t-\sigma}^\t\,\phi\,dt}.
\eeq
Thus, the inverse problem to identify $\g$ from the exact observation $\phi$ satisfying \eqref{exact_obs_data} has a unique solution $\g\in L^\infty(\partial\O)$ given in \eqref{gamma_rep}. 

Next, we observe that the inverse problem {\bf (IP)} is non linear since the temperature distribution $u(x,t)$ depends on the impedance coefficient $\g$. Also, it is ill-posed in the sense that small perturbation in the observation data $\phi$ in \eqref{exact_obs_data} may lead to large deviation in the corresponding solution of the inverse problem. This is the case because the expression \eqref{gamma_rep} for $\g$ contains a derivative of the observation $\phi$. Thus, in order to obtain some stable approximations for $\g$ some regularization method has to be employed. We will consider the Tikhonov regularization for obtaining stable approximations in the next section. 

\section{Convergence rates with weak source condition}\label{sec-convergence}
Let $u(\g)$ denotes the unique weak solution of \eqref{pde1} for a fixed $\g\in \mathcal{A}$. For $\d>0$, let $\phi_\d$ be the noisy data corresponding to the exact data $\phi$ satisfying \eqref{noisy_obs_error}. As discussed in the previous section, the inverse problem is ill-posed, and thus we shall use the Tikhonov-regularization in order to obtain stable approximations. Throughout we shall denote by $\g^\dagger\in L^2(\partial\O)$ the exact impedance coefficient to be identified for the corresponding exact data $\phi.$ For a fixed $\a>0$, consider the output-least square Tikhonov functional
\beq\label{Tikh_funct}
J(\g):=\int_{\t-\sigma}^\t\int_{\O}|u(\g)-\phi_{\d}|^2\,dx\,dt+\a \|\g-\g^*\|_{L^2(\partial\O)}^2,
\eeq
where $\g^*\in L^2(\partial\O)$ is an initial guess for $\g^\dagger$ that incorporates some apriori smoothness assumption on $\g^\dagger.$ By $\g^\a_\d$ we denote a minimizer of the optimization problem
\beq\label{Tikh_minimization}
\min_{\g\in\mathcal{A}} J(\g).
\eeq
These minimizers are the Tikhonov-regularized solutions. Note that the initial guess $\g^*$ may not belong to the admissible set $\mathcal{A}$. Our next two results are about the existence and stability of such minimizers, which in turn ensures that the minimizers are indeed regularized solutions. It is to be noted that by now the arguments for the proof of existence and stability is well established in the literature (cf. \cite{engl_zou_2000}), but since the context of this paper is different, in order to keep the paper self contained we include the proof also.
\bt\label{Tikh_minimizer_existence}
The minimization problem \eqref{Tikh_minimization} has a solution.
\et
\bpf
It is clear that $\mathcal{A}$ is a convex set. Let $\{\g_n\}$ be a minimizing sequence in $\mathcal{A}.$ Clearly, $\{\g_n\}$ is a bounded sequence in $L^2(\partial\O).$ Thus, there exists a subsequence $\{\g_m\}$ and a $\g^\a\in L^2(\partial\O)$ such that $\g_m\rightharpoonup \g^\a$ in $L^2(\partial\O).$ Now the closedness and convexity of $\mathcal{A}$ implies that $\mathcal{A}$ is weakly closed, thus, $\g^\a\in \mathcal{A}.$ Since $u(\g_m)$ is  the weak solution of \eqref{pde1} for $\g=\g_m,$ by the estimate \eqref{weak_sol_estimate}, it follows that $\{u(\g_m)\}$ is a bounded sequence in $L^2(0,\t;H^1(\O)).$ Thus, there exists a subsequence, still denoted by $\{u(\g_m)\}$, and $u^*\in L^2(0,\t;H^1(\O))$ such that $u(\g_m)\rightharpoonup u^*$ in $L^2(0,\t;H^1(\O))$ as $m\to \infty.$ Now, for any $\eta\in L^2(0,\t;H^1(\O))$, writing 
$$\int_0^\t\int_{\partial\O}\g_m\,u(\g_m)\eta\,dx\,dt=\int_0^\t\int_{\partial\O}\g^\a\, u(\g_m)\eta\,dx\,dt+\int_0^\t\int_{\partial\O}(\g_m-\g^\a)\,u(\g_m)\eta\,dx\,dt,$$
using the weak convergence of $\g_m$ and the boundedness of $u(\g_m)$, it follows that $$\int_0^\t\int_{\partial\O}\g_m\,u(\g_m)\eta\,dx\,dt\to\int_0^\t\int_{\partial\O}\g^\a\,u^*\eta\,dx\,dt\q\text{as}\,\,m\to \infty.$$
Therefore, using the fact that $u(\g_m)$ is a weak solution of \eqref{pde1} for $\g=\g_m$, and the weak convergence of $\{u(\g_m)\}$, it follows that 
\beqarray
\int_0^\t\int_\O[-u^*\frac{\partial\eta}{\partial t}+\nabla u^*\cdot\nabla \eta]\,dx\,dt+\int_0^\t\int_{\partial\O}\g^\a\,u^*\eta\,dx\,dt=\int_0^\t\int_\O f\eta\,dx\,dt+\int_0^\t\int_{\partial\O}g\eta\,dx\,dt+\int_{\O}u_0\eta(\cdot,0)\,dx
\eeqarray
for all $\eta\in H^1(0,\t;H^1(\O))$ with $\eta(\cdot,\t)=0.$ Since $u(\g^\a)$ is the unique weak solution of \eqref{pde1} for $\g=\g^\a$, we have $u^*=u(\g^\a).$
We now consider the identity
\beqarray
&&\liminf_m\int_{\t-\sigma}^\t\int_{\O}|u(\g_m)-\phi_\d|^2\,dx\,dt\\
&&\q =\liminf_m\int_{\t-\sigma}^\t\int_{\O}\Big[|u(\g_m)-u(\g^\a)|^2+|u(\g^\a)-\phi_\d|^2+2\big(u(\g_m)-u(\g^\a)\big)\big(u(\g^\a)-\phi_\d\big)\Big]\,dx\,dt.
\eeqarray
Since $u(\g_m)$ and $u(\g^\a)$ are the weak solutions of \eqref{pde1} for $\g=\g_m$ and $\g^\a$, respectively, using the estimate \eqref{weak_sol_estimate}, we have
$$\lim_{m\to\infty}\int_{\t-\sigma}^\t\int_{\O}|u(\g_m)-u(\g^\a)|^2\,dx\,dt=0.$$
Therefore, $$\liminf_m\int_{\t-\sigma}^\t\int_{\O}|u(\g_m)-\phi_\d|^2\,dx\,dt=\int_{\t-\sigma}^\t\int_{\O}|u(\g^\a)-\phi_\d|^2\,dx\,dt.$$
Thus, using the weak lower semi-continuity of the $L^2$-norm and using the fact that $\{\g_n\}$ is a minimizing sequence for the minimization problem \eqref{Tikh_minimization}, we have
\beqarray
&&\int_{\t-\sigma}^\t\int_{\O}|u(\g^\a)-\phi_\d|^2\,dx\,dt+\a\|\g^\a-\g^*\|^2_{L^2(\partial\O)}\\
&&\q\leq \liminf_m\Big\{\int_{\t-\sigma}^\t\int_{\O}|u(\g_m)-\phi_\d|^2\,dx\,dt+\a\|\g_m-\g^*\|^2_{L^2(\partial\O)}\Big\}\\
&&\q = \min_{\g\in \mathcal{A}}J(\g).
\eeqarray
This shows that $\g^\a$ is a minimizer of \eqref{Tikh_minimization}.
\epf
We now prove the stability of the minimization problem \eqref{Tikh_minimization} with respect to the observation data $\phi_\d.$ That is, we establish that the minimizers of $\eqref{Tikh_minimization}$ are indeed regularized solutions. 
\bt\label{Tikh_funct_stability}
Let $\{\phi_n\}$ be a sequence such that $\phi_n$ converges to $\phi_\d$ in $L^2(\t-\sigma,\t;L^2(\O)).$ For a fixed $\a>0$, let $\g_n^\a$ be the minimizer of 
\beq\label{pert_Tikh_funct}
\min_{\g\in\mathcal{A}}\int_{\t-\sigma}^\t\int_{\O}|u(\g)-\phi_n|^2\,dx\,dt+\a\|\g-\g^*\|^2_{L^2(\partial\O)}.
\eeq
Then there exists a subsequence $\{\g_n^\a\}$ that converges to a minimizer $\g_\d^\a.$
\et
\bpf
Since $\g_n^\a$ is a minimizer of \eqref{pert_Tikh_funct}, we have
$$
\int_{\t-\sigma}^\t\int_{\O}|u(\g_n^\a)-\phi_n|^2\,dx\,dt+\a\|\g_n^\a-\g^*\|^2_{L^2(\partial\O)}\leq \int_{\t-\sigma}^\t\int_{\O}|u(\g)-\phi_n|^2\,dx\,dt+\a\|\g-\g^*\|^2_{L^2(\partial\O)}
$$
for any $\g\in \mathcal{A}.$ Thus, $\{\g_n^\a\}$ is a bounded sequence in $L^2(\partial\O).$ Therefore, there exists a subsequence, with abuse of notation, denoted by $\{\g_n^\a\}$ and a $\g^\a\in L^2(\partial\O)$ such that $\g_n^\a$ converges weakly to $\g^\a$ in $L^2(\partial\O)$ as $n\to\infty.$ Now, the closedness and convexity of $\mathcal{A}$ implies that $\g^\a\in \mathcal{A}.$ 

Since $u(\g_n^\a)$ is the unique weak solution of \eqref{pde1}, by \eqref{weak_sol_estimate}  it follows that $\{u(\g_n^\a)\}$ is bounded in $L^2(0,\t;H^1(\O)).$ Thus, there exist a subsequence, still denoted by $u(\g_n^\a)$ and a $u_*\in L^2(0,\t;H^1(\O))$ such that $u(\g_n^\a)$ converges weakly to $u_*$ in $L^2(0,\t;H^1(\O))$ as $n\to \infty.$ 

Now, for any $\eta\in L^2(0,\t;H^1(\O))$, writing 
$$\int_0^\t\int_{\partial\O}\g_n^\a\,u(\g_n^\a)\eta\,dx\,dt=\int_0^\t\int_{\partial\O}\g^\a\, u(\g_n^\a)\eta\,dx\,dt+\int_0^\t\int_{\partial\O}(\g_n^\a-\g^\a)\,u(\g_n^\a)\eta\,dx\,dt,$$
using the weak convergence of $\g_n^\a$ and the boundedness of $u(\g_n^\a)$, it follows that $$\int_0^\t\int_{\partial\O}\g_n^\a\,u(\g_n^\a)\eta\,dx\,dt\to\int_0^\t\int_{\partial\O}\g^\a\,u_*\eta\,dx\,dt\q\text{as}\,\,n\to \infty.$$
Therefore, using the fact that $u(\g_n^\a)$ is a weak solution, and the weak convergence of $\{u(\g_n^\a)\}$, it follows that 
\beqarray
&&\int_0^\t\int_\O[-u_*\frac{\partial\eta}{\partial t}+\nabla u_*\cdot\nabla \eta]\,dx\,dt+\int_0^\t\int_{\partial\O}\g^\a\,u_*\eta\,dx\,dt\\
&&\q=\int_0^\t\int_\O f\eta\,dx\,dt+\int_0^\t\int_{\partial\O}g\eta\,dx\,dt+\int_{\O}u_0\eta(\cdot,0)\,dx
\eeqarray
for all $\eta\in H^1(0,\t;H^1(\O))$ with $\eta(\cdot,\t)=0.$ Since $u(\g^\a)$ is the unique weak solution of \eqref{pde1}, we have $u_*=u(\g^\a).$

Now the weak lower semi continuity of $L^2$-norm implies
$$\|\g^\a-\g^*\|^2_{L^2(\partial\O)}\leq \liminf_n\,\|\g_n^\a-\g^*\|^2_{L^2(\partial\O)}$$
and
$$\int_{\t-\sigma}^\t\int_{\O}|u(\g^\a)-\phi_\d|^2\,dx\,dt\leq \liminf_n\int_{\t-\sigma}^\t\int_{\O}|u(\g_n^\a)-\phi_\d|^2\,dx\,dt.$$
Since $\phi_n\to \phi_\d$ in $L^2(\t-\sigma,\t;L^2(\O)),$ it follows that
\beqarray
&&\liminf_n \int_{\t-\sigma}^\t\int_{\O}|u(\g_n^\a)-\phi_n|^2\,dx\,dt \\
&&\q =\liminf_n \int_{\t-\sigma}^\t\int_{\O}\Big[|u(\g_n^\a)-\phi_\d|^2+|\phi_n-\phi_\d|^2+2\big(u(\g_n^\a)-\phi_\d\big)(\phi_\d-\phi_n)\Big]\,dx\,dt\\
&&\q =\liminf_n\int_{\t-\sigma}^\t\int_{\O}|u(\g_n^\a)-\phi_\d|^2\,dx\,dt\\
&&\q =\liminf_n\int_{\t-\sigma}^\t\int_{\O}\Big[|u(\g_n^\a)-u(\g^\a)|^2+|u(\g^\a)-\phi_\d|^2+2\big(u(\g_n^\a)-u(\g^\a)\big)\big(u(\g^\a)-\phi_\d\big)\Big]\,dx\,dt.
\eeqarray
Since $u(\g_n^\a)$ and $u(\g^\a)$ are the weak solutions of \eqref{pde1} for $\g=\g_n^\a$ and $\g^\a$, respectively, using the estimate \eqref{weak_sol_estimate}, we have
$$\lim_{n\to\infty}\int_{\t-\sigma}^\t\int_{\O}|u(\g_n^\a)-u(\g^\a)|^2\,dx\,dt=0.$$
Therefore, $$\liminf_n\int_{\t-\sigma}^\t\int_{\O}|u(\g_n^\a)-\phi_n|^2\,dx\,dt=\int_{\t-\sigma}^\t\int_{\O}|u(\g^\a)-\phi_\d|^2\,dx\,dt.$$
Thus,
\beqarray
&&\int_{\t-\sigma}^\t\int_{\O}|u(\g^\a)-\phi_\d|^2\,dx\,dt+\a\|\g^\a-\g^*\|^2_{L^2(\partial\O)}\\
&&\q\leq \liminf_n\Big\{\int_{\t-\sigma}^\t\int_{\O}|u(\g_n^\a)-\phi_n|^2\,dx\,dt+\a\|\g_n^\a-\g^*\|^2_{L^2(\partial\O)}\Big\}\\
&&\q\leq \limsup_n\Big\{\int_{\t-\sigma}^\t\int_{\O}|u(\g_n^\a)-\phi_n|^2\,dx\,dt+\a\|\g_n^\a-\g^*\|^2_{L^2(\partial\O)}\Big\}\\
&&\q\leq \limsup_n\Big\{\int_{\t-\sigma}^\t\int_{\O}|u(\g)-\phi_n|^2\,dx\,dt+\a\|\g-\g^*\|^2_{L^2(\partial\O)}\Big\}\\
&&\q=\int_{\t-\sigma}^\t\int_{\O}|u(\g)-\phi_\d|^2\,dx\,dt+\a\|\g-\g^*\|^2_{L^2(\partial\O)}
\eeqarray
for any $\g\in \mathcal{A}.$ This shows that $\g^\a$ is a minimizer, and we denote this by $\g_\d^\a.$ Also, taking $\g=\g_\d^\a$ in the last equality, we have
\beq\label{sm}
\begin{cases}
&\lim_n\int_{\t-\sigma}^\t\int_{\O}|u(\g_n^\a)-\phi_n|^2\,dx\,dt+\a\|\g_n^\a-\g^*\|^2_{L^2(\partial\O)}\\
&\q=\int_{\t-\sigma}^\t\int_{\O}|u(\g^\a_\d)-\phi_\d|^2\,dx\,dt+\a\|\g^\a_\d-\g^*\|^2_{L^2(\partial\O)}
\end{cases}
\eeq

We now establish the convergence of $\g_n^\a$ to $\g_\d^\a$ in $L^2(\partial\O)$ as $n\to \infty.$ The proof is by contradiction. Suppose that $\g_n^\a$ does not converge to $\g_\d^\a$ in $L^2(\partial\O).$ Then, clearly
$$\|\g_\d^\a-\g^*\|^2_{L^2(\partial\O)}<\limsup_n\|\g_n^\a-\g^*\|^2_{L^2(\partial\O)}=:\epsilon.$$
Thus, there exists a subsequence of $\{\g_n^\a\}$, say $\{\g_m^\a\}$ such that $\g_m^\a$ converges weakly to $\g_\d^\a$ in $L^2(\partial\O)$ and $$\|\g_m^\a-\g^*\|^2_{L^2(\partial\O)}\to \epsilon\q\text{as}\,\,m\to \infty.$$
Now from \eqref{sm}, we have
\beqarray
&&\lim_m\int_{\t-\sigma}^\t\int_{\O}|u(\g_m^\a)-\phi_m|^2\,dx\,dt\\
&&\q=\int_{\t-\sigma}^\t\int_{\O}|u(\g^\a_\d)-\phi_\d|^2\,dx\,dt+\a\Big(\|\g_\d^\a-\g^*\|^2_{L^2(\partial\O)}-\lim_m\|\g_m^\a-\g^*\|^2_{L^2(\partial\O)}\Big)\\
&&\q=\int_{\t-\sigma}^\t\int_{\O}|u(\g^\a_\d)-\phi_\d|^2\,dx\,dt +\a\underbrace{\Big(\|\g_\d^\a-\g^*\|^2_{L^2(\partial\O)}-\epsilon\Big)}_{< 0}.
\eeqarray
Hence,
$$\lim_m\int_{\t-\sigma}^\t\int_{\O}|u(\g_m^\a)-\phi_m|^2\,dx\,dt<\int_{\t-\sigma}^\t\int_{\O}|u(\g^\a_\d)-\phi_\d|^2\,dx\,dt.$$
Therefore, we have
\beqarray
\int_{\t-\sigma}^\t\int_{\O}|u(\g^\a_\d)-\phi_\d|^2\,dx\,dt &=& \liminf_n\int_{\t-\sigma}^\t\int_{\O}|u(\g_n^\a)-\phi_n|^2\,dx\,dt\\
&\leq & \lim_m\int_{\t-\sigma}^\t\int_{\O}|u(\g_m^\a)-\phi_m|^2\,dx\,dt \\
&<&\int_{\t-\sigma}^\t\int_{\O}|u(\g^\a_\d)-\phi_\d|^2\,dx\,dt,
\eeqarray
 a contradiction.
\epf
Let $f\in L^2(0,\t;L^2(\O)), g\in L^2(0,\t;L^2(\partial\O))$ and $u_0\in L^2(\O)$. Then by Theorem \ref{existence_estimate_weak_sol}, we know that for each $\g\in \mathcal{A}$, \eqref{pde1} has a unique weak solution in $L^2(0,\t;H^1(\O))\cap L^\infty(0,\t;L^2(\O))$ and we have used the notation $u(\g)$ to denote this weak solution. For notational simplicity we denote this by a map $F$, that is, $F:\mathcal{A}\to L^2(0,\t;H^1(\O))\cap L^\infty(0,\t;L^2(\O))$ is defined as $\g\mapsto F(\g):=u(\g).$ 

In order to avoid notational complexity, throughout $C$ will denote a generic constant that may depend only on $\O,\t,\sigma, \underline{\g}, \overline{\g}, u_0,f, g,$ the operator norm of the trace operator (and also on some function $\psi$ considered later in Theorem \ref{convergence_rate}).
\bl\label{frechet_diff}
The mapping $F:\mathcal{A}\subset L^\infty(\partial\O)\to L^2(0,\t;H^1(\O))\cap L^\infty(0,\t;L^2(\O))$ is Fr{\'e}chet differentiable.
\el
\bpf
Since $u(\g)$ denotes the unique weak solution of \eqref{pde1}, we have
\beq\label{weak_frecehet}
\begin{cases}
&\int_{0}^\tau \int_\O  [-u(\g)\frac{\partial\eta}{\partial t}+\nabla u(\g)\cdot\nabla \eta]\,dx\,dt+\int_{0}^\tau\int_{\partial\O}\g u(\g)\eta\,dx\,dt\\
&\q=\int_{0}^\tau \int_\O f\eta\,dx\,dt+\int_{0}^\tau\int_{\partial\O}g\eta\,dx\,dt+\int_\O\,u_0\eta(\cdot,0)\,dx
\end{cases}
\eeq
for all $\eta\in H^1(0,\tau;H^1(\O))$ with $\eta(\cdot,\tau)=0.$
Let $h\in L^\infty(\partial\O)$ be such that $\g+h\in \mathcal{A}.$ 
 Then, we have
\beq\label{weak_frecehet_2}
\begin{cases}
&\int_{0}^\tau \int_\O  [-u(\g+h)\frac{\partial\eta}{\partial t}+\nabla u(\g+h)\cdot\nabla \eta]\,dx\,dt+\int_{0}^\tau\int_{\partial\O}\g u(\g+h)\eta\,dx\,dt\\
&\q=\int_{0}^\tau \int_\O f\eta\,dx\,dt+\int_{0}^\tau\int_{\partial\O}g\eta\,dx\,dt+\int_\O\,u_0\eta(\cdot,0)\,dx
\end{cases}
\eeq
for all $\eta\in H^1(0,\tau;H^1(\O))$ with $\eta(\cdot,\tau)=0.$ Thus, from \eqref{weak_frecehet} and \eqref{weak_frecehet_2}, we have
\beq\label{weak_frecehet_3}
\begin{cases}
&\int_{0}^\tau \int_\O  -[u(\g+h)-u(\g)]\frac{\partial\eta}{\partial t}\,dx\,dt+\int_{0}^\t\int_\O[\nabla u(\g+h)-\nabla u(\g)]\cdot\nabla \eta\,dx\,dt\\
&\q +\int_{0}^\tau\int_{\partial\O}(\g+h)[ u(\g+h)-u(\g)]\eta\,dx\,dt+ \int_{0}^\t\int_{\partial\O}\,hu(\g)\eta\,dx\,dt=0
\end{cases}
\eeq
for all $\eta\in H^1(0,\tau;H^1(\O))$ with $\eta(\cdot,\tau)=0.$
Now consider the PDE
\beq\label{weak_frechet_auxilary_PDE}
\begin{cases}
\frac{\partial v}{\partial t}=\Delta v \q &  \text{in}\q\O_\t,\\
\frac{\partial v}{\partial \nu}+\g v =-hu(\g)\q & \text{on}\q\partial\O_\t,\\
v(\cdot,0)=0\q & \text{in}\q\O.
\end{cases}
\eeq
By Theorem \ref{existence_estimate_weak_sol} there exist a unique weak solution $v$ of \eqref{weak_frechet_auxilary_PDE} satisfying
\beq\label{frechet_auxiliary_integral_1}
\int_{0}^\t\int_\O v\frac{\partial\eta}{\partial t}\,dx\,dt-\int_{0}^\t\int_\O \nabla v\cdot \nabla\eta\,dx\,dt-\int_{0}^\t\int_{\partial\O}\g v \eta\,dx\,dt =\int_{0}^\t\int_{\partial\O} hu(\g)\eta\,dx\,dt
\eeq
for all $\eta\in h^1(0;\t,H^1(\O))$ with $\eta(\cdot,\t)=0$ in $\O.$ Moreover, from \eqref{weak_frechet_auxilary_PDE} and the estimate \eqref{weak_sol_estimate} it follows that 
\beqarray
\max_{t\in[0,\t]}\|v(t)\|_{L^2(\O)}+\|v\|_{L^2(0,\t;H^1(\O))}&\leq & C\|hu(\g)\|_{L^2(0,\t;L^2(\partial\O))}\\
&\leq & C\|h\|_{L^\infty(\partial\O)}\|u(\g)\|_{L^2(0,\t;L^2(\partial\O))}\\
&\leq & C\|h\|_{L^\infty(\partial\O)}\|u(\g)\|_{L^2(0,\t;H^1(\O))}.
\eeqarray 
Therefore, for a fixed $\g\in \mathcal{A},$ the map $h\mapsto v,$ where $v$ is the unique weak solution of \eqref{weak_frechet_auxilary_PDE}, is a bounded linear operator from $L^\infty(\partial\O)\to L^\infty(0,\t;L^2(\O))\cap L^2(0,\t;H^1(\O)).$
Let $w:=u(\g+h)-u(\g)-v$, then from \eqref{frechet_auxiliary_integral_1} and \eqref{weak_frecehet_3}, we have
\beq\label{frechet_auxiliary_integral_2}
\int_{0}^\t\int_\O -w\frac{\partial\eta}{\partial t}\,dx\,dt+\int_{0}^\t\int_{\partial\O}(\g+h) w \eta\,dx\,dt +\int_{0}^\t\int_{\partial\O} hv\eta\,dx\,dt
\eeq
for all $\eta\in H^1(0,\t;H^1(\O))$ with $\eta(\cdot,\t)=0$ in $\O.$ Therefore from \eqref{frechet_auxiliary_integral_2}, the definition of weak solution and the estimate \eqref{weak_sol_estimate}, it follows that there exists a constant $C$ depending only on $\underline{\g},\,\O$ and $\t$ such that
\beq\label{frechet_auxilary_estimate}
\max_{t\in [0,\t]}\|w(t)\|_{L^2(\O)}+\|w\|_{L^2(0,\t;H^1(\O))}\leq C\|hv\|_{L^2(0,\t;L^2(\partial\O))}.
\eeq

Now using the fact that $u(\g)$ is a weak solution of \eqref{pde1}, from Theorem \ref{existence_estimate_weak_sol} it follows that
$$\|u(\g)\|_{L^2(0,\t;H^1(\O))}\leq  C\|h\|_{L^\infty(\partial\O)}(\|f\|_{L^2(0,\t;L^2(\O))}+\|g\|_{L^2(0,\t;L^2(\partial\O))}+\|u_0\|_{L^2(\O)}).$$
Therefore, from \eqref{frechet_auxilary_estimate} we have
$$
\begin{cases}
&\max_{t\in [0,\t]}\|w(t)\|_{L^2(\O)}+\|w\|_{L^2(0,\t;H^1(\O))}\\
&\q\leq C\|h\|^2_{L^\infty(\partial\O)}(\|f\|_{L^2(0,\t;L^2(\O))}+\|g\|_{L^2(0,\t;L^2(\partial\O))}+\|u_0\|_{L^2(\O)}).
\end{cases}
$$
Thus, it follows that $F$ is Fr{\'e}chet differentiable, and $$F'(\g)h=v,$$ where $v$ is the unique weak solution of \eqref{weak_frechet_auxilary_PDE} for a given $\g\in \mathcal{A}.$
\epf
We are now in a position to state and proof our main result about the convergence rate under a weak source condition. Let us recall that $\g^\dagger\in \mathcal{A}$ is the exact Robin coefficient to be identified and $\g^*\in L^2(\partial\O)$ is an apriori initial guess for $\g^\dagger$, and $u(\g^\dagger)$ denotes the unique weak solution of \eqref{pde1} for $\g=\g^\dagger.$
\bt{\rm{\bf (Convergence rate)}}\label{convergence_rate}
Let $\psi\in H^1_0(\t-\sigma,\t;L^2(\O))\cap\, L^2(\t-\sigma,\t; H^2(\O))$ be such that 
\beq\label{source_cond}
\int_{\t-\sigma}^\t\,u(\g^\dagger)\psi\,dt=\g^\dagger\,-\g^*\q\text{on}\,\,\partial\O.
\eeq
For fixed $\a,\d>0,$ let $\g_\d^\a$ be the minimizer of \eqref{Tikh_funct}.
Then $$\|\g^\dagger-\g^\a_\d\|_{L^2(\partial\O)}\leq C\Big(\frac{\d^2}{\a}+\a\Big)^\frac{1}{2},$$
for some constant $C>0$ independent of $\a,\,\,\d, \g^\dagger$ and $\g^\a_\d.$
\et
\bpf
Since $\g_\d^\a$ is a minimizer of \eqref{Tikh_funct} for a fixed $\a,$ using \eqref{noisy_obs_error}, it follows that
\beqarray
&&\int_{\t-\sigma}^\t\int_{\O}|u({\g_\d^\a})-\phi_{\d}|^2\,dx\,dt+\a \|\g_\d^\a-\g^*\|_{L^2(\partial\O)}\\
&&\q\leq \int_{\t-\sigma}^\t\int_{\O}|u({\g^\dagger})-\phi_{\d}|^2\,dx\,dt+\a \|\g^\dagger-\g^*\|_{L^2(\partial\O)}\\
&&\q \leq \d^2+\a\|\g^\dagger-\g^*\|^2_{L^2(\O)}.
\eeqarray
Thus,
\beqarray
&&\int_{\t-\sigma}^\t\|u({\g_\d^\a})-\phi_\d\|^2_{L^2(\O)}+\a\|\g^\dagger-\g_\d^\a\|^2_{L^2(\partial\O)}\\
&&\q\leq \d^2+\a\big[\|\g^\dagger-\g^\a_\d\|^2_{l^2(\partial\O)}+\|\g^\dagger-\g^*\|^2_{L^2(\partial\O)}-\|\g^\a_\d-\g^*\|^2_{L^2(\partial\O)}\big]\\
&&\q =\d^2 +2\a\<\g^\dagger-\g^\a_\d, \g^\dagger-\g^*\>_{L^2(\partial\O)}.
\eeqarray
The source condition implies
$$\a\<\g^\dagger-\g^\a_\d,\g^\dagger-\g^*\>_{L^2(\partial\O)}=\a\<\int_{\t-\sigma}^\t u({\g^\dagger})\psi\,dt,\g^\dagger-\g^\a_\d\>_{L^2(\partial\O)}=\a\int_{\t-\sigma}^\t\int_{\partial\O}u({\g^\dagger})\psi\, (\g^\dagger-\g^\a_\d)\,dx\,dt.$$
Since $F$ is Fr{\'e}chet differentiable, taking $h=\g^\a_\d-\g^\dagger,$ $\eta=\psi,$ $u(\g)=u(\g^\dagger)$ and using the identity \eqref{frechet_auxiliary_integral_1}, we get
\beqarray
\int_{\t-\sigma}^\t\int_{\partial\O}u(\g^\dagger)\psi\, (\g^\a_\d-\g^\dagger)\,dx\,dt &=&\int_{\t-\sigma}^\t\int_\O\,[F'(\g^\dagger)(\g^\a_\d-\g^\dagger)\frac{\partial\psi}{\partial t}-\nabla F'(\g^\dagger)(\g^\a_\d-\g^\dagger)\cdot\nabla \psi]\,dx\,dt\\
&&\q -\int_{\t-\sigma}^\t\int_{\partial\O}\g^\dagger F'(\g^\dagger)(\g^\a_\d-\g^\dagger)\psi \,dx\,dt.
\eeqarray
Thus,
\beqarray
\a\<\g^\dagger-\g^\a_\d,\g^\dagger-\g^*\>_{L^2(\partial\O)}&=&\a\Big(\int_{\t-\sigma}^\t\int_\O\,\big[-F'(\g^\dagger)(\g^\a_\d-\g^\dagger)\frac{\partial\psi}{\partial t}+\nabla F'(\g^\dagger)(\g^\a_\d-\g^\dagger)\cdot\nabla \psi\big]\,dx\,dt\\
&&\q\q+\int_{\t-\sigma}^\t\int_{\partial\O}\g^\dagger F'(\g^\dagger)(\g^\a_\d-\g^\dagger)\psi \,dx\,dt\Big).
\eeqarray
Let $w_\d^\a:=u(\g^\a_\d)-u(\g^\dagger)-F'(\g^\dagger)(\g^\a_\d-\g^\dagger).$ Then, we have
\beq\label{auxilary_reduced_source}
\begin{cases}
&\a\<\g^\dagger-\g^\a_\d,\g^\dagger-\g^*\>_{L^2(\partial\O)}\\
&\q=\a\Big[\int_{\t-\sigma}^\t\int_\O\,\big(-[u(\g^\a_\d)-u(\g^\dagger)]\frac{\partial\psi}{\partial t}+\nabla[u(\g^\a_\d)-u(\g^\dagger)]\cdot\nabla\psi\big)\,dx\,dt \\
&\q\q\q+\int_{\t-\sigma}^\t\int_{\partial\O}\g^\dagger[u(\g^\a_\d)-u(\g^\dagger)]\psi\,dx\,dt\\
&\q\q\q+\int_{\t-\sigma}^\t\int_{\O}\big(w_\d^\a\frac{\partial\psi}{\partial t}-\nabla w^\a_\d\cdot \nabla \psi\big)\,dx\,dt-\int_{\t-\sigma}^\t\int_{\partial\O}\g^\dagger w^\a_\d \psi\,dx\,dt\Big]
\end{cases}
\eeq
Taking $w=w^\a_\d$, from \eqref{frechet_auxiliary_integral_2}, we have
\beqarray
\int_{\t-\sigma}^\t\int_{\O}[-w_\d^\a\frac{\partial\psi}{\partial t}+\nabla w_\d^\a\cdot \nabla\psi]\,dx\,dt+\int_{\t-\sigma}^\t\int_{\partial\O}\g^\a_\d w^\a_\d\psi\,dx\,dt+\int_{\t-\sigma}^\t\int_{\partial\O}(\g^\a_\d-\g^\dagger)F'(\g^\dagger)(\g^\a_\d-\g^\dagger)\psi\,dx\,dt=0.
\eeqarray
Thus,
\beqarray
\int_{\t-\sigma}^\t\int_{\O}[w_\d^\a\frac{\partial\psi}{\partial t}-\nabla w_\d^\a\cdot \nabla\psi]\,dx\,dt-\int_{\t-\sigma}^\t\int_{\partial\O}\g^\dagger w^\a_\d\psi\,dx\,dt=\int_{\t-\sigma}^\t\int_{\partial\O}(\g^\a_\d-\g^\dagger)(u(\g^\a_\d)-u(\g^\dagger))\psi\,dx\,dt.
\eeqarray
Therefore, from \eqref{auxilary_reduced_source}, we have
\beqarray
&&\a\<\g^\dagger-\g^\a_\d,\g^\dagger-\g^*\>_{L^2(\partial\O)}\\
&&\q =\a\Big[\int_{\t-\sigma}^\t\int_\O\,\big(-[u(\g^\a_\d)-u(\g^\dagger)]\frac{\partial\psi}{\partial t}+\nabla[u(\g^\a_\d)-u(\g^\dagger)]\cdot\nabla\psi\big)\,dx\,dt\\
&&\q\q\q+\int_{\t-\sigma}^\t\int_{\partial\O}\g^\a_\d(u(\g^\a_\d)-u(\g^\dagger))\psi\,dx\,dt\Big]\\
&&\q = \a\Big[\int_{\t-\sigma}^\t\int_\O\,\big(-[u(\g^\a_\d)-u(\g^\dagger)]\frac{\partial\psi}{\partial t}+[u(\g^\dagger)-u(\g^\a_\d)]\Delta\psi\big)\,dx\,dt\\
&&\q\q\q+\int_{\t-\sigma}^\t\int_{\partial\O}[u(\g^\a_\d)-u(\g^\dagger)]\frac{\partial\psi}{\partial\nu}\,dx\,dt+\int_{\t-\sigma}^\t\int_{\partial\O}\g^\a_\d[u(\g^\a_\d)-u(\g^\dagger)]\psi\,dx\,dt\Big]\\
&&\q=: I_1+I_2+I_3+I_4.
\eeqarray
We now estimate the integrals $I_1, I_2, I_3$ and $I_4.$ We will be using the Cauchy-Schwarz inequality and Young's inequality appropriately which will involve an arbitrary $\e>0.$
\beqarray
|I_1|&=&|\a\int_{\t-\sigma}^\t\int_\O\,-[u(\g^\a_\d)-u(\g^\dagger)]\frac{\partial\psi}{\partial t}\,dx\,dt|\\
&\leq &\a\int_{\t-\sigma}^\t\|u(\g^\a_\d)-u(\g^\dagger)\|_{L^2(\O)}\|\frac{\partial\psi}{\partial t}\|_{L^2(\O)}\,dt\\
&\leq & \a\int_{\t-\sigma}^\t\|u(\g^\a_\d)-\phi_\d\|_{L^2(\O)}\|\frac{\partial\psi}{\partial t}\|_{L^2(\O)}\,dt+\a\int_{\t-\sigma}^\t\|\phi_\d-u(\g^\dagger)\|_{L^2(\O)}\|\frac{\partial\psi}{\partial t}\|_{L^2(\O)}\,dt\\
&\leq & \e\int_{\t-\sigma}^\t\|u(\g^\a_\d)-\phi_\d\|^2_{L^2(\O)}\,dt+\frac{\a^2}{4\e}\int_{\t-\sigma}^\t\|\frac{\partial\psi}{\partial t}\|^2_{L^2(\O)}\,dt \\
&&\q+ \e\int_{\t-\sigma}^\t\|\phi_\d-u(\g^\dagger)\|^2_{L^2(\O)}+\frac{\a^2}{4\e}\int_{\t-\sigma}^\t\|\frac{\partial\psi}{\partial t}\|^2_{L^2(\O)}\,dt\\
&\leq &\e\int_{\t-\sigma}^\t\|u(\g^\a_\d)-\phi_\d\|^2_{L^2(\O)}\,dt+\d^2\e+\frac{\a^2}{2\e}\int_{\t-\sigma}^\t\|\frac{\partial\psi}{\partial t}\|^2_{L^2(\O)}\,dt.
\eeqarray
Similarly, we obtain
\beqarray
|I_2|&=& |\a\int_{\t-\sigma}^\t[u(\g^\dagger)-u(\g^\a_\d)]\Delta\psi\,dx\,dt|\\
&\leq & \e\int_{\t-\sigma}^\t\|u(\g^\a_\d)-\phi_\d\|^2_{L^2(\O)}\,dt+\d^2\e+\frac{\a^2}{2\e}\int_{\t-\sigma}^\t\|\Delta\psi\|^2_{L^2(\O)}\,dt.
\eeqarray
In order to estimate $I_3$ and $I_4$, additionally we will make use of the continuity of the trace map from $H^1(\O)$ to $L^2(\partial\O).$
\beqarray
|I_3|&=& |\a\int_{\t-\sigma}^\t\int_{\partial\O}[u(\g^\a_\d)-u(\g^\dagger)]\frac{\partial\psi}{\partial\nu}\,dx\,dt|\\
&\leq & \a\int_{\t-\sigma}^\t\|u(\g^\a_\d)-u(\g^\dagger)\|_{L^2(\partial\O)}\|\frac{\partial\psi}{\partial\nu}\|_{L^2(\partial\O)}\,dt\\
&\leq &C\a\int_{\t-\sigma}^\t\|u(\g_\d^\a)-u(\g^\dagger)\|_{L^2(\O)}\|\frac{\partial\psi}{\partial\nu}\|_{L^2(\partial\O)}\,dt\\
&\leq &C\a\int_{\t-\sigma}^\t\|u(\g^\a_\d)-\phi_\d\|_{L^2(\O)}\|\frac{\partial\psi}{\partial \nu}\|_{L^2(\partial\O)}\,dt+C\a\int_{\t-\sigma}^\t\|\phi_\d-u(\g^\dagger)\|_{L^2(\O)}\|\frac{\partial\psi}{\partial \nu}\|_{L^2(\partial\O)}\,dt\\
&\leq & \e\int_{\t-\sigma}^\t\|u(\g_\d^\a)-u(\g^\dagger)\|^2_{L^2(\O)}\,dt+\d^2\e+\frac{C^2\a^2}{2\e}\int_{\t-\sigma}^\t\|\frac{\partial\psi}{\partial\nu}\|^2_{L^2(\partial\O)}\,dt.
\eeqarray
Similarly, we have
\beqarray
|I_4|&=&|\a\int_{\t-\sigma}^\t\int_{\partial\O}\g^\a_\d[u(\g_\d^\a)-u(\g^\dagger)]\psi\,dx\,dt|\\
&\leq & \overline{\g}\,C\a\int_{\t-\sigma}^\t\|u(\g^\a_\d)-u(\g^\dagger)\|_{L^2(\O)}\|\psi\|_{L^2(\partial\O)}\,dx\,dt\\
&\leq & \e\int_{\t-\sigma}^\t\|u(\g^\a_\d)-\phi_\d\|^2_{L^2(\O)}\,dt+\d^2\e+\frac{\overline{\g}^2\,C^2\a^2}{2\e}\int_{\t-\sigma}^\t\|\psi\|^2_{L^2(\partial\O)}\,dt.
\eeqarray
Therefore, we obtain
\beqarray
&&\int_{\t-\sigma}^\t\|u(\g_\d^\a)-\phi_\d\|^2_{L^2(\O)}+\a\|\g^\dagger-\g^\a_\d\|^2_{L^2(\partial\O)}\\
&&\q\leq \d^2+8\e\int_{\t-\sigma}^\t\|u(\g_\d^\a)-\phi_\d\|^2_{L^2(\O)}\,dt+8\d^2\e\\
&&\q\q+\frac{\a^2}{\e}\Big[\int_{\t-\sigma}^\t\Big(\|\frac{\partial\psi}{\partial t}\|^2_{L^2(\O)}+\|\Delta\psi\|^2_{L^2(\O)}+C^2\|\frac{\partial\psi}{\partial\nu}\|^2_{L^2(\partial\O)}+\overline{\g}^2C^2\|\psi\|^2_{L^2(\partial\O)}\Big)\,dt\Big].
\eeqarray
Thus, taking $\e=\frac{1}{16},$ we have
\beqarray
&&\int_{\t-\sigma}^\t\|u(\g_\d^\a)-\phi_\d\|^2_{L^2(\O)}\,dt+2\a\|\g^\dagger-\g^\a_\d\|^2_{L^2(\partial\O)}\\
&&\q\leq 3\d^2+32\a^2\Big[\int_{\t-\sigma}^\t\Big(\|\frac{\partial\psi}{\partial t}\|^2_{L^2(\O)}+\|\Delta\psi\|^2_{L^2(\O)}+C^2\|\frac{\partial\psi}{\partial\nu}\|^2_{L^2(\partial\O)}+\overline{\g}^2C^2\|\psi\|^2_{L^2(\partial\O)}\Big)\,dt\Big]\\
&&\q\leq C(\d^2+\a^2).
\eeqarray
Therefore,
$$\int_{\t-\sigma}^\t\|u(\g_\d^\a)-\phi_\d\|^2_{L^2(\O)}\,dt\leq C(\d^2+\a^2)$$
and
$$\|\g^\dagger-\g_\d^\a\|_{L^2(\partial\O)}\leq C\Big(\frac{\d^2}{\a}+\a\Big)^{\frac{1}{2}}.$$
\epf
\brem
From the estimate obtained in the above theorem for $\|\g^\dagger-\g^\a_\d\|_{L^2(\partial\O)}$ under the source condition \eqref{source_cond} it follows that if we choose the regularization parameter $\a$ as $\a\sim\d$, then we have $\|\g^\dagger-\g_\d^\a\|_{L^2(\partial\O)}=O(\d^\frac{1}{2}).$
\erem
\brem
As observed, our analysis could produce the rate $O(\d^\frac{1}{2})$. It would be interesting to analyse the convergence rates by proceeding along the recent line of research based on variational source conditions (see for e.g., \cite{chen_yousept_2021, hohage_weidling_2015, hohage_weidling_2017_siam, sprung_hohage_2019}) which to the best of our knowledge still remains to be explored in the context of impedance coefficient identification. In fact we are aware of only a very recent paper \cite{chen_jiang_yousept_zou_2022} that deals with Robin coefficient identification using variational source condition for elliptic PDE only. 
\erem

\section{Discussion about source condition}\label{sec-source_diss}
We now discuss advantages of the source condition \eqref{source_cond}. From Theorem \ref{convergence_rate} it is clear that our source condition does not require any smallness condition to be verified unlike the smallness condition of the form \eqref{source_smallness},  that is required in standard convergence theory for Tikhonov-regularization for non linear operators (cf. \cite{engl_hanke_neubauer, engl_kunisch_neubauer_1989}). Moreover, in order to obtain a convergence rate of $O(\d^\frac{1}{2})$ we do not need to assume that $\g^\dagger$ belongs to the range of the adjoint of the Fr{\'e}chet derivative of the parameter-to-solution operator. This is another big advantage since most often the range of such operators is nothing but some Sobolev spaces with higher smoothness (see for e.g. \cite{hohage_1997}), and thus our source condition being free from such apriori range condition means that we can obtain the said convergence rate without imposing higher smoothness assumption on $\g^\dagger,$ which is obviously not known to us. 

We now discuss about the apriori regularity of $\g^\dagger$ that is embedded in the source condition \eqref{source_cond}. Recall that for defining trace operator, the minimal assumption that we need is that $\int_{\t-\sigma}^\t\,u(\g^\dagger)\,\psi\,dt\in H^1(\O)$, and in that case it follows that we must have $\g^\dagger-\g^*\in H^\frac{1}{2}(\partial\O)$. Thus, in order that the source condition \eqref{source_cond} makes sense, the regularity assumption that we require is that $\g^\dagger-\g^*\in H^\frac{1}{2}(\partial\O).$ We shall show that in our setting we do always have $\int_{\t-\sigma}^\t\,u(\g^\dagger)\psi\,dt\in H^1(\O)$, thanks to the following very recent result on multiplication of elements in Sobolev spaces \cite{behzadan_holst_2021}. 

\bt\label{pointwise-multiplication} (cf. \cite[Theorem 7.4]{behzadan_holst_2021})
Let $D\subset\R^d,\,d\in\{1,2,3\},$ be a bounded domain with Lipschitz continuous boundary. Let $s_i(i=1,2),\,s$ be real numbers such that $s_i\geq s\geq 0$ and $s_1+s_2-s>\frac{d}{2}$. Then there exists a constant $C$ depending only on $s_1,s_2,s,d$ and $\O$ such that $$\|UV\|_{H^s(D)}\leq C\,\|U\|_{H^{s_1}(D)}\|V\|_{H^{s_2}(D)},\q\forall\,U\in H^{s_1}(D),\,\,V\in H^{s_2}(D).$$
\et
Since $\psi\in L^2(\t-\sigma,\t;H^2(\O))$ and $u(\g^\dagger)\in L^2(\t-\sigma,\t;H^1(\O))$, by Theorem \ref{pointwise-multiplication} it follows that $\int_{\t-\sigma}^\t\,u(\g^\dagger)\psi\,dt\in H^1(\O)$. 

\brem
It is to be noted that one can also derive a sufficient condition for $\int_{\t-\sigma}^\t\,u(\g^\dagger)\psi\,dt\in H^1(\O)$ based on the results related to Banach algebra properties of Sobolev space, but with this argument we need to assume $H^2$- spatial regularity of $u(\g^\dagger)$, which of course is a drawback. Indeed, recall that for $d\in \{2,3\}$, $\O\subset\R^d$, $H^s(\O),\,s\geq 2,$ is a Banach algebra. Thus, if we assume that $u(\g^\dagger)\in L^2(\t-\sigma,\t;H^2(\O))$(such a regularity indeed holds if the initial profile $u_0$ and the source function $f$ in \eqref{pde1} are assumed to be of appropriate spatial regularity) then $\int_{\t-\sigma}^\t\,u(\g^\dagger)\,\psi\,dt\in H^2(\O)\subset H^1(\O)$ and thus by the property of trace operator it follows that we must have $\g^\dagger-\g^*\in H^\frac{3}{2}(\partial\O)\subset H^\frac{1}{2}(\partial\O).$
\erem

We now look into the compatibility issue related to the source condition. From the expression of the source condition \eqref{source_cond}, one may think that on some portion of the boundary $\partial\O$ near the terminal time status, it may happen that $u(\g^\dagger)$ vanishes(in the sense of trace) but $\g^\dagger-\g^*$ does not, and this implies that on such portions one has to know $\g^\dagger.$ But this is compatible with the fact that for terminal time status, if $u(\g^\dagger)$ vanishes (in the sense of trace) on some portion of the boundary $\partial\O$ then it is impossible to recover $\g^\dagger$ on such portions, as can be observed from \eqref{gamma_rep}. 

Next, we give a procedure to verify the source condition explicitly, motivated from the construction given in \cite{engl_zou_2000} for the case of $1$-dimensional diffusion coefficient identification problem. However, our construction is valid for the dimension $d\in\{2,3\}$ as considered in this paper, but provided we have $u(\g^\dagger)\in L^2(\t-\sigma,\t;H^2(\O))$. Indeed, we construct a $\psi\in H^1_0(\t-\sigma,\t;L^2(\O))\cap L^2(\t-\sigma,\t;H^2(\O))$ that satisfies the source condition \eqref{source_cond}. Before proceeding further, we recall the following interesting result from \cite{sprung_thesis}.
\bt(cf. \cite[Theorem A.2.9]{sprung_thesis})\label{reciprocal}
Let $D\subset\R^d$ be a bounded domain with Lipschitz boundary. Let $s>\frac{d}{2}$. If $U\in H^s(D)$ and $U\geq c>0$ for some constant $c,$ then $\frac{1}{U}\in H^s(D).$
\et
Our next result is about verification of source condition under certain smoothness assumption.
\bt\label{source_verify_construct}
Let $\g^\dagger,\g^*\in L^\infty(\partial\O)\cap H^\frac{3}{2}(\partial\O)$ and $u(\g^\dagger)\in L^2(\t-\sigma,\t;H^2(\O))$. Let $\psi_1\in H^1(\t-\sigma,\t)$ be an arbitrary function. Define $$u_1(x):=\int_{\t-\sigma}^\t(\t-t)(\t-\sigma-t)\psi_1\,u(\g^\dagger)\,dt.$$ Assume that 
\beq\label{assumption_source_verify}
u_1\neq 0\,\, \text{a.e. in}\,\, \O\q\text{and}\q \frac{1}{u_1}\in L^\infty(\O). 
\eeq
Then there exists a $\psi\in H^1_0(\t-\sigma,\t;L^2(\O))\cap L^2(\t-\sigma,\t;H^2(\O))$ satisfying the source condition \eqref{source_cond}.
\et
\bpf
Let $\psi_2\in H^1(\t-\sigma,\t)$ and $\f\in H^m(\O)$, $m\geq 2,$ be two arbitrary functions. We will write $u$ to denote $u(\g^\dagger).$ Define $$u_2(x):=\int_{\t-\sigma}^\t(\t-t)(\t-\sigma-t)\psi_2\f\,u\,dt.$$
Let $\Upsilon^\dagger, \Upsilon^*\in H^2(\O)$ be such that $\g^\dagger$ and $\g^*$ are their traces, respectively. We define
$$\psi:=(\t-t)(\t-\sigma-t)\Big(\psi_2\f+\psi_1\frac{\Upsilon^\dagger-\Upsilon^*-u_2}{u_1}\Big).$$
We claim that this $\psi$ satisfies the source condition. Observe that from the $L^\infty$-assumption in \eqref{assumption_source_verify} and from Theorem \ref{reciprocal}, it follows that $\frac{1}{u_1}\in H^2(\O)$ and hence $\frac{\Upsilon^\dagger-\Upsilon^*-u_2}{u_1}\in L^2(\t-\sigma,\t;H^2(\O)).$ Also, from the definition of $\psi$ it clearly follows that $\psi(\cdot,\t-\sigma)=0=\psi(\cdot,\t)$ in $\O$. Thus, $\psi\in H^1_0(\t-\sigma,\t;L^2(\O))\cap L^2(\t-\sigma,\t;H^2(\O)).$ Moreover,
\beqarray
\int_{\t-\sigma}^\t\,u(\g^\dagger)\psi\,dt&=&\int_{\t-\sigma}^\t(\t-t)(\t-\sigma-t)\Big(\psi_2\f\,u+\psi_1\,u\frac{\g^\dagger-\g^*-u_2}{u_1}\Big)\q\text{ on}\,\,\partial\O\\
&=& \g^\dagger-\g^*\q\text{on}\,\,\partial\O.
\eeqarray
This completes the proof.
\epf
\brem
 Note that the non-zero assumption in \eqref{assumption_source_verify} does make sense, because otherwise it would mean $u(\g^\dagger)$ vanishes on the boundary near the terminal time $\t$(i.e., $[\t-\sigma,\t]$), but then from \eqref{gamma_rep} it follows that the impedance coefficient $\g^\dagger$ is impossible to recover.
\erem
\section{Conclusion}
We have considered an inverse problem of identifying a spatially dependent impedance coefficient in a parabolic PDE from a short-time observation of the temperature distribution. We have proposed a weak-type source condition that allowed to obtain a convergence rate of $O(\d^\frac{1}{2})$ provided we choose the regularization parameter $\a\sim\d$, where $\d$ is the deterministic noise level. As compared to the standard convergence theory of Tikhonov regularization for non linear operators in Hilbert spaces, the above rate of convergence is obtained under a simple source condition which does not require any range condition of the adjoint of Fr{\'e}chet derivative, no smallness condition to be verified and requires a regularity assumption that $\g^\dagger\in H^\frac{1}{2}(\partial\O).$ Moreover, if it is known that $u(\g^\dagger)$ vanishes(in the sense of trace) on some portion of the boundary $\partial\O,$ then one has to know $\g^\dagger$ apriori over those portions, as it is impossible to recover the impedance coefficient on those portions.

\noi
{\bf Acknowledgments.} The author is grateful to Prof. Thorsten Hohage (University of G{\"o}ttingen) for several insightful comments and for bringing the result in \cite{sprung_thesis} to our notice, which has helped to improve the content. This work was initiated when the author was a postdoctoral researcher in the Institute for Numerical and Applied Mathematics, University of G{\"o}ttingen. The present financial support from the TIFR Centre for Applicable Mathematics, as a post doctoral fellow is gratefully acknowledged.


\begin{thebibliography}{99}
\bibitem{adams} R.\,Adams, {\it Sobolev spaces.} Pur Appl Math. 65, Academic press, New York, 1975.

\bibitem{bellassoued_cheng_choulli_2008} M.\,Bellassoued, J.\,Cheng and M.\,Choulli, {\it Stability estimate for an inverse boundary coefficient problem in thermal imaging.} J. Math. Anal. Appl. 343 (2008), no. 1, 328–336.

\bibitem{beck_blackwell_clair} J.V.\,Beck, B.\,Blackwell and C.R.S.\,Clair, {\it Inverse Heat Conduction: Ill-Posed Problems,} Wiley-Interscience, New York, 1985.

\bibitem{behzadan_holst_2021}  A.\,Behzadan and M.\,Holst,  {\it Multiplication in Sobolev spaces, revisited.} Ark. Mat. 59 (2021), no. 2, 275–306. 



\bibitem{cao_pereverzev_2006} H.\,Cao and S.V.\,Pereverzev, {\it Natural linearization for the identification of a diffusion coefficient in a quasi-linear parabolic system from short-time observations.} Inverse Problems 22 (2006), 2311-2330.

\bibitem{k.cao_2022} K.\,Cao, {\it Convergence rates for the reaction coefficient and the initial temperature identification problems.} Appl. Anal. 101 (2022), no. 7, 2472–2497.

\bibitem{chen_yousept_2021} D.H.\,Chen and I.\,Yousept, {\it Variational source conditions in $L^p$-spaces.} SIAM J. Math. Anal. 53 (2021), no. 3, 2863–2889. 

\bibitem{chen_jiang_yousept_zou_2022} D.H.\,Chen, D.\,Jiang, I.\,Yousept and J.\,Zou, {\it Variational source conditions for inverse Robin and flux problems by partial measurements.} Inverse Problems and Imaging, 2022, 16(2): 283-304. 

\bibitem{cheng_lu_yamamoto_2012} J.\,Cheng, S.\,Lu and M.\,Yamamoto, {\it Reconstruction of the Stefan-Boltzmann coefficients in a heat-transfer process.} Inverse Problems 28 (2012), no. 4, 045007, 17 pp.



\bibitem{engl_hanke_neubauer} H.W.\,Engl, M.\,Hanke and A.\,Neubauer, {\it Regularization of inverse problems.} Kluwer, Dordrecht, 1996.

\bibitem{engl_kunisch_neubauer_1989} H.W.\,Engl, K.\,Kunisch and A.\,Neubauer, {\it Convergence rates for Tikhonov regularisation of nonlinear ill-posed problems.} Inverse Problems 5 (1989), no. 4, 523–540.

\bibitem{engl_zou_2000}  H.W.\,Engl and J.\,Zou, {\it A new approach to convergence rate analysis of Tikhonov regularization for parameter identification in heat conduction}, Inverse Problems 16 (2000) 1907-1923.


\bibitem{evans}  L.C.\,Evans, {\it Partial differential equations.} Second edition, Graduate Studies in Mathematics, 19. American Mathematical Society, Providence, RI, 2010. xxii+749 pp. ISBN: 978-0-8218-4974-3 35-01.



\bibitem{hao_thanh_lesnic_2013} D.N.\, H{\`a}o, P.X.\,Thanh and D.\,Lesnic, {\it Determination of the heat transfer coefficients in transient heat conduction.} Inverse Problems 29 (2013), no. 9, 095020, 21 pp.

\bibitem{hohage_1997} T.\,Hohage, {\it Logarithmic convergence rates of the iteratively regularized Gauss-Newton method for an inverse potential and an inverse scattering problem.} Inverse Problems 13 (1997), no. 5, 1279–1299.

\bibitem{hohage_weidling_2015} T.\,Hohage and F.\,Weidling, {\it Verification of a variational source condition for acoustic inverse medium scattering problems.} Inverse Problems 31 (2015), no. 7, 075006, 14 pp.

\bibitem{hohage_weidling_2017_siam} T.\,Hohage and F.\,Weidling, {\it Characterizations of variational source conditions, converse results, and maxisets of spectral regularization methods.} SIAM J. Numer. Anal. 55 (2017), no. 2, 598–620.


\bibitem{isakov_1991} V.\,Isakov, {\it Inverse Parabolic Problems with the final overdetermination.} Comm. Pure Appl. Math., 44 (1991), 185-209.

\bibitem{jin_lu_2012} B.\,Jin and X.\,Lu, {\it Numerical identification of a Robin coefficient in parabolic problems.} Math. Comp. 81 (2012), no. 279, 1369–1398.

\bibitem{jin_zou_2009} B.\,Jin and J.\,Zou, {\it Numerical estimation of piecewise constant Robin coefficient.} SIAM J. Control Optim. 48 (2009), no. 3, 1977–2002.



\bibitem{keung_zou_1998} Y.L.\,Keung and J.\,Zou, {\it Numerical identifications of parameters in parabolic systems.} Inverse Problems 14 (1998), no. 1, 83–100.

\bibitem{kugler_sincich_2009} P.\,K{\"u}gler and E.\,Sincich, {\it Logarithmic convergence rates for the identification of a nonlinear Robin coefficient.} J. Math. Anal. Appl. 359 (2009), no. 2, 451–463. 

\bibitem{liu_wang_2016} J.\,Liu and Y.\,Wang, {\it On the reconstruction of boundary impedance of a heat conduction system from nonlocal measurement.} Inverse Problems 32 (2016), no. 7, 075002, 23 pp.



\bibitem{nairopeq} M.T. Nair, {\it Linear Operator Equations: Approximation and Regularization}. World Scientific, Hackensack, 2009.


\bibitem{neubauer_1989} A.\,Neubauer, {\it Tikhonov regularisation for nonlinear ill-posed problems: optimal convergence rates and finite-dimensional approximation.} Inverse Problems 5 (1989), no. 4, 541–557.



\bibitem{schock_1984} E.\,Schock, {\it On the asymptotic order of accuracy of Tikhonov regularization.} J. Optim. Theory Appl. 44 (1984), no. 1, 95–104.

\bibitem{sprung_thesis} B.\,Sprung, {\it Convergence rates for variational regularization of statistical inverse problems.} Ph.D. Thesis, University of G{\"o}ttingen (2019).

\bibitem{sprung_hohage_2019} B.\,Sprung and T.\,Hohage, {\it Higher order convergence rates for Bregman iterated variational regularization of inverse problems.} Numer. Math. 141 (2019), no. 1, 215–252.

\bibitem{troltzsch} F.\,Tr{\"o}ltzsch, {\it Optimal control of partial differential equations. Theory, methods and applications.} Graduate Studies in Mathematics, 112. American Mathematical Society, Providence, RI, 2010. xvi+399 pp. ISBN: 978-0-8218-4904-0

\bibitem{wang_liu_2017} Y.\,Wang and J.\,Liu, {\it On the simultaneous recovery of boundary heat transfer coefficient and initial heat status.}  J. Inverse Ill-Posed Probl. 25 (2017), no. 5, 597–616.

\bibitem{yamamoto_zou_2001} M.\,Yamamoto and J.\,Zou, {\it Simultaneous reconstruction of the initial temperature and heat radiative coefficient. } Inverse Problems 17 (2001) 1181–1202.

\bibitem{zhang_liu_2021} M.\,Zhang and J.\,Liu, {\it On the simultaneous reconstruction of boundary Robin coefficient and internal source in a slow diffusion system.} Inverse Problems 37 (2021), no. 7, Paper No. 075008, 33 pp.


\end{thebibliography}
\end{document}